\documentclass{amsart}
\usepackage{hyperref}
\newcommand*{\mailto}[1]{\href{mailto:#1}{\nolinkurl{#1}}}

\newtheorem{theorem}{Theorem}[section]

\newtheorem{lemma}[theorem]{Lemma}
\newtheorem{corollary}[theorem]{Corollary}
\newtheorem{remark}[theorem]{Remark}
\newtheorem{hypothesis}[theorem]{Hypothesis}

\newcommand{\R}{{\mathbb R}}
\newcommand{\N}{{\mathbb N}}
\newcommand{\Z}{{\mathbb Z}}
\newcommand{\C}{{\mathbb C}}

\newcommand{\nn}{\nonumber}
\newcommand{\be}{\begin{equation}}
\newcommand{\ee}{\end{equation}}

\newcommand{\ol}{\overline}
\newcommand{\ti}{\tilde}

\newcommand{\id}{{\rm 1\hspace{-0.6ex}l}}
\newcommand{\E}{\mathrm{e}}
\newcommand{\I}{\mathrm{i}}

\newcommand{\im}{\mathrm{Im}}
\newcommand{\re}{\mathrm{Re}}
\newcommand{\dom}{\mathfrak{D}}

\newcommand{\floor}[1]{\lfloor#1 \rfloor}
\newcommand{\ceil}[1]{\lceil#1 \rceil}

\newcommand{\eps}{\varepsilon}

\newcommand{\sig}{\sigma}
\newcommand{\lam}{\lambda}
\newcommand{\gam}{\gamma}
\newcommand{\om}{\omega}

\newcommand{\Lam}{\Lambda}
\newcommand{\ta}{\theta}
\newcommand{\wt}{\widetilde}

\newcommand{\cA}{\mathcal{A}}

\newcommand{\Om}{\Omega}
\newcommand{\pa}{\partial}

\numberwithin{equation}{section}

\begin{document}

\title[Commutation Methods for Strongly Singular Potentials]{Commutation Methods for Schr\"odinger Operators with Strongly Singular Potentials}

\author[A.\ Kostenko]{Aleksey Kostenko}
\address{Institute of Applied Mathematics and Mechanics\\
NAS of Ukraine\\ R. Luxemburg str. 74\\
Donetsk 83114\\ Ukraine\\ and School of Mathematical Sciences\\
Dublin Institute of Technology\\
Kevin Street\\ Dublin 8\\ Ireland}
\email{\mailto{duzer80@gmail.com}}

\author[A.\ Sakhnovich]{Alexander Sakhnovich}
\address{Faculty of Mathematics\\ University of Vienna\\
Nordbergstrasse 15\\ 1090 Wien\\ Austria}
\email{\mailto{Oleksandr.Sakhnovych@univie.ac.at}}
\urladdr{\url{http://www.mat.univie.ac.at/~sakhnov/}}

\author[G.\ Teschl]{Gerald Teschl}
\address{Faculty of Mathematics\\ University of Vienna\\
Nordbergstrasse 15\\ 1090 Wien\\ Austria\\ and International
Erwin Schr\"odinger
Institute for Mathematical Physics\\ Boltzmanngasse 9\\ 1090 Wien\\ Austria}
\email{\mailto{Gerald.Teschl@univie.ac.at}}
\urladdr{\url{http://www.mat.univie.ac.at/~gerald/}}

\thanks{Math. Nachr. (to appear)}
\thanks{{\it Research supported by the Austrian Science Fund (FWF) under Grant No.\ Y330}}

\keywords{Schr\"odinger operators, spectral theory, commutation methods, strongly singular potentials}
\subjclass[2010]{Primary 34B20, 34L05; Secondary 34B24, 47A10}

\begin{abstract}
We explore the connections between singular Weyl--Titchmarsh theory and the single and
double commutation methods. In particular, we compute the singular Weyl function of the
commuted operators in terms of the original operator. We apply the results to spherical
Schr\"odinger operators (also known as Bessel operators). We also investigate the connections
with the generalized B\"acklund--Darboux transformation.
\end{abstract}

\maketitle

\section{Introduction}
\label{sec:int}

The present paper is concerned with spectral theory for one-dimensional Schr\"o\-dinger operators
\be
H = -\frac{d^2}{dx^2} + q(x), \qquad x\in (a,b),
\ee
on the Hilbert space $L^2(a,b)$ with a real-valued potential $q\in L^1_{\mathrm{loc}}(a,b)$. It has been
shown recently by Gesztesy and Zinchenko \cite{gz}, Fulton and Langer \cite{ful08}, \cite{fl}, Kurasov and Luger \cite{kl}
that, for a large class of singularities at $a$, it is still possible to define a singular Weyl function at the basepoint $a$.
Furthermore, in previous work we have shown that this singular Weyl function shares many properties with the classical
 Weyl function \cite{kst2} and established the connection with super singular perturbations for the
 special case of spherical Schr\"odinger operators
 \be
H = -\frac{d^2}{dx^2} + \frac{l(l+1)}{x^2} + q(x), \qquad x\in (0,\infty),
\ee
(also known as Bessel operators) \cite{kt}.

On the other hand, commutation methods have played an important role in the
theory of one-dimensional Schr\"o\-dinger operators both as a method for inserting eigenvalues as
well as for constructing solutions of the (modified) Korteweg--de Vries equation (see, e.g., \cite{gss}
and the references therein).
Historically, these methods of inserting eigenvalues go back to Jacobi
\cite{ja} and Darboux \cite{da} with decisive later contributions by Crum
\cite{cr}, Krein \cite{kr57}, Schmincke \cite{sc}, and Deift \cite{de}. Two
particular methods turned out to be of special importance: The
single commutation method, also called the Crum--Darboux method \cite{cr},
\cite{da} (actually going back at least to Jacobi \cite{ja}) and the double
commutation method, to be found, e.g., in the seminal work of Gel'fand and
Levitan \cite{gl}. For recent extensions of these methods we refer to
\cite{fg}, \cite{gstdef}, \cite{gt}, \cite{sc2}.

Krein \cite{kr57} was the first to realize the connection between inverse spectral problems and Crum's results.
Namely, in \cite{kr57}, the connection between the spectral measures of the original and transformed operators
was established and then exploited to characterize the spectral measures of Bessel operators in the case
$l\in\N$ (see also \cite{fad}). This idea has been subsequently used by many authors: See, for instance,
\cite{ahm2, car, de, fad, gks} and references therein.

The purpose of our present paper is to continue the work of Krein and establish the connection between
the singular Weyl functions (and hence between the spectral measures) of the original and transformed
operators for both the single and double commutation method. In particular, we will obtain an independent
proof for the fact that the singular Weyl function of perturbed  Bessel operators is a generalized Nevanlinna function.
In addition, we investigate the connections with the generalized B\"acklund--Darboux transformation (GBDT)
for a particular example. This method is a generalization of the double commutation method which it contains
as a special case (cf.\ Subsection~\ref{Comp}).

\section{Singular Weyl--Titchmash theory}

We begin by recalling a few facts from \cite{kst2}.
To set the stage, we will consider one-dimensional Schr\"odinger operators on $L^2(a,b)$
with $-\infty \le a<b \le \infty$ of the form
\begin{equation} \label{stli}
\tau = - \frac{d^2}{dx^2} + q,
\end{equation}
where the potential $q$ is real-valued and satisfies
\begin{equation}
q \in L^1_{loc}(a,b).
\end{equation}
We will use $\tau$ to denote the formal differential expression and $H$ to denote a corresponding
self-adjoint operator given by $\tau$ with separated boundary conditions at $a$ and/or $b$.

If $a$ (resp.\ $b$) is finite and $q$ is in addition integrable
near $a$ (resp.\ $b$), we will say $a$ (resp.\ $b$) is a \textit{regular}
endpoint.  We will say $\tau$, respectively $H$, is \textit{regular} if
both $a$ and $b$ are regular.

We will choose a point $c\in(a,b)$ and also consider the operators $H^D_{(a,c)}$, $H^D_{(c,b)}$
which are obtained by restricting $H$ to $(a,c)$, $(c,b)$ with a Dirichlet boundary condition at
$c$, respectively. The corresponding operators with a Neumann boundary condition will be
denoted by $H^N_{(a,c)}$ and $H^N_{(c,b)}$.

Moreover, let $c(z,x)$, $s(z,x)$ be the solutions of $\tau u = z\, u$ corresponding
to the initial conditions $c(z,c)=1$, $c'(z,c)=0$ and $s(z,c)=0$, $s'(z,c)=1$.

Define the Weyl $m$-functions (corresponding to the base point $c$) such that
\begin{align}\nn
u_-(z,x) &= c(z,x) - m_-(z) s(z,x), \qquad z\in\C\setminus\sig(H^D_{(a,c)}),\\\label{defupm}
u_+(z,x) &= c(z,x) + m_+(z) s(z,x), \qquad z\in\C\setminus\sig(H^D_{(c,b)}),
\end{align}
are square integrable on $(a,c)$, $(c,b)$ and satisfy the boundary condition of $H$
at $a$, $b$ (if any), respectively. The solutions $u_\pm(z,x)$ (as well as their multiples)
are called Weyl solutions at $a$, $b$.
For further background we refer to \cite[Chap.~9]{tschroe} or \cite{wdln}.

To define an analogous singular Weyl $m$-function at the, in general singular, endpoint $a$ we will
first need an analog of the system of solutions $c(z,x)$ and $s(z,x)$. Hence our first goal is
to find a system of real entire solutions $\theta(z,x)$ and $\phi(z,x)$ such that $\phi(z,x)$ lies in the domain of
$H$ near $a$ and such that the Wronskian $W(\theta(z),\phi(z))=1$. By a real entire function we mean
an entire function which is real-valued on the real line. To this end we start with a
hypothesis which will turn out necessary and sufficient for such a system of solutions to exist.

\begin{hypothesis}\label{hyp:gen}
Suppose that the spectrum of $H^D_{(a,c)}$ is purely discrete for one (and hence for all)
$c\in (a,b)$.
\end{hypothesis}

Note that this hypothesis is for example satisfied if $q(x) \to +\infty$ as $x\to a$ (cf.\ Problem~9.7 in \cite{tschroe}).

\begin{lemma}[\cite{kst2}]\label{lem:pt}
Suppose Hypothesis~\ref{hyp:gen} holds. Then there exists a fundamental system of solutions $\phi(z,x)$ and
$\theta(z,x)$ of $\tau u = z u$ which are real entire with respect to $z$ such that
\be\label{eq:wronsk}
W(\theta(z),\phi(z)) = 1
\ee
and $\phi(z,.)$ is in the domain of $H$ near $a$.
Here $W_x(u,v)= u(x)v'(x)-u'(x)v(x)$ is the usual Wronski determinant.
\end{lemma}

It is important to point out that such a system is not unique and any other such sytem is given by
\[
\ti{\theta}(z,x) = \E^{-g(z)} \theta(z,x) - f(z) \phi(z,x), \qquad
\ti{\phi}(z,x) = \E^{g(z)} \phi(z,x),
\]
where $g(z)$, $f(z)$ are real entire functions.

Given a system of real entire solutions $\phi(z,x)$ and $\theta(z,x)$ as in the above lemma we can define the
singular Weyl function
\be\label{defM}
M(z) = -\frac{W(\theta(z),u_+(z))}{W(\phi(z),u_+(z))}
\ee
such that the solution which is in the domain of $H$ near $b$ (cf.\ \eqref{defupm}) is given by
\be
u_+(z,x)= a(z) \big(\theta(z,x) + M(z) \phi(z,x)\big),
\ee
where $a(z)= - W(\phi(z),u_+(z)) = \beta(z) - m_+(z) \alpha(z)$.
By construction we obtain that the singular Weyl function $M(z)$ is analytic in $\C\backslash\R$ and satisfies $M(z)=M(z^*)^*$.
Rather than $u_+(z,x)$ we will use
\be\label{defpsi}
\psi(z,x)= \theta(z,x) + M(z) \phi(z,x).
\ee

Recall also from \cite[Lem.~3.2]{kst2} that associated with $M(z)$ is a corresponding spectral measure
\be\label{defrho}
\frac{1}{2} \left( \rho\big((x_0,x_1)\big) + \rho\big([x_0,x_1]\big) \right)=
\lim_{\eps\downarrow 0} \frac{1}{\pi} \int_{x_0}^{x_1} \im\big(M(x+\I\eps)\big) dx.
\ee

\section{Connection with the single commutation method}

\subsection{Preliminary basic results}
We begin by recalling a few basic facts from the single commutation method.
Let $A$ be a densely defined closed operator and recall that $H =A^* A$ is a self-adjoint
operator with $\ker(H)=\ker(A)$. Similarly, $\hat{H}=A A^*$ is a self-adjoint
operator with $\ker(\hat{H})=\ker(A^*)$. Then the key observation is the following
well-known result (see e.g.\ \cite[Thm.~8.6]{tschroe} for a short proof):

\begin{theorem}[\cite{de}]
Let $A$ be a densely defined closed operator and introduce $H= A^* A$, $\hat{H}= A A^*$.
Then the operators $H\big|_{\ker(H)^\perp}$ and $\hat{H}\big|_{\ker(\hat{H})^\perp}$ are unitarily
equivalent.

If $H \psi =  \lam \psi$, $\lam\in\R$, $\psi\in\dom(H)$, then $\hat{\psi}=A \psi \in\dom(\hat{H})$ with
$\hat{H} \hat{\psi} = \lam \hat{\psi}$ and $\|\hat{\psi}\|=\sqrt{\lam} \|\psi\|$. Moreover,
\be\label{resasa}
R_{\hat{H}}(z) \supseteq \frac{1}{z} \left(A R_H(z) A^* - \id\right), \quad
R_H(z) \supseteq \frac{1}{z} \left(A ^* R_{\hat{H}}(z) A - \id\right),
\ee
where $R_H(z)=(H-z)^{-1}$ denotes the resolvent of an operator $H$.
\end{theorem}

\subsection{Application to Schr\"odinger operators}
In order to find such a factorization for a given Schr\"odinger operator $H$ one requires a positive solution of the underlying differential
equation. Existence of such a solution is equivalent to semi-boundedness of $H$ and we will thus make the following assumption:

\begin{hypothesis}\label{hyp:sc}
In addition to Hypothesis~\ref{hyp:gen} assume that $H$ is bounded from below and limit point at $b$. Let $\lam\in\R$ be such that
$H-\lam \ge 0$.
\end{hypothesis}

In particular, let $\phi(z,x)$, $\theta(z,x)$ be a fundamental system of solutions as in Lemma~\ref{lem:pt} and recall (cf.\ \cite[Lem.~9.7]{tschroe}) that the
Green's function of $H$ is given by
\be
G(z,x,y)= \begin{cases}
\phi(z,x) \psi(z,y), & x \le y,\\
\phi(z,y) \psi(z,x), & x \ge y.
\end{cases}
\ee
Moreover, if $H - \lam \ge 0$ we must have $\phi(\lam,x)>0$ (as well as $\psi(\lam,x)>0$) for $x\in(a,b)$ possibly after flipping signs by \cite[Cor.~2.4]{gstz}.

Now consider the operator
\begin{align}\nn
A_\phi f &= a_\phi f, \qquad a_\phi = - \frac{d}{dx} + \frac{\phi'(\lam,x)}{\phi(\lam,x)}, \\
\dom(A_\phi) &= \{ f\in L^2(a,b) | f\in AC_{loc}(a,b),\: a_\phi f \in L^2(a,b) \}.
\end{align}
Here we use $A_\phi$ and $a_\phi$ for the operator and differential expression, respectively.
It is straightforward to check (cf.\ \cite[Problem~9.3]{tschroe}) that $A_\phi$ is closed and that its adjoint
is given by
\begin{align}\nn
A_\phi^* f &=  a^*_\phi f, \qquad a^*_\phi = \frac{d}{dx} + \frac{\phi'(\lam,x)}{\phi(\lam,x)}, \\ \nn
\dom(A_\phi^*) &= \{ f\in L^2(a,b) | f\in AC_{loc}(a,b),\: a^*_\phi f \in L^2(a,b),\\
&\qquad \lim_{x\to a,b} f(x) g(x) =0, \forall g\in\dom(A_\phi) \}.
\end{align}
If we also have $\theta(\lam,x)>0$ we can also define $A_\theta$ by using $\theta(\lam,x)$ in place
of $\phi(\lam,x)$.

\begin{lemma}
Assume Hypothesis~\ref{hyp:sc} holds. Then $H - \lam = A_\phi^* A_\phi$. If in addition $\theta(\lam,x)>0$ and
$\tau$ is limit point at $a$, then we also have $H - \lam = A_\theta^* A_\theta$.
\end{lemma}

\begin{proof}
It is simple algebra to check that $H$ and $A_\phi^* A_\phi$ agree on functions with compact support in $(a,b)$ hence $A_\phi^* A_\phi$ is
a self-adjoint extension of the minimal operator associated with $\tau$ and it remains to identify the boundary conditions.
Since by assumption $\tau$ is limit point at $b$ we only need to consider $a$. Moreover, since $\phi(\lam,x)$ is in the domain
of $A_\phi^* A_\phi$ near $a$ by our choice of $A_\phi$, both $A_\phi^* A_\phi$ and $H$ are associated with the boundary condition generated by
$\phi(\lam,x)$ if $\tau$ is limit circle at $a$.

The case of $A_\theta$ is even simpler since by our limit point assumption there is no need to identify any boundary conditions.
\end{proof}

The commuted operator $\hat{H} - \lam = A_\phi A_\phi^*$ is associated with the potential
\be\label{eq:hat_q}
\hat{q}(x) = q(x) - 2 \frac{d}{dx} \frac{\phi'(\lam,x)}{\phi(\lam,x)}.
\ee
Moreover, it is straightforward to check that if $u(z)$ solves $\tau u = z u$ then
\be\label{eq:hatu}
\hat{u}(z,x) = a_\phi u(z,x) = -\frac{W_x(\phi(\lam),u(z))}{\phi(\lam,x)}
\ee
solves $\hat{\tau} \hat{u} = z \hat{u}$ and given two solutions $u(z)$ and $v(z)$ we have
\be\label{eq:wronski}
W(\hat{u}(z),\hat{v}(z)) = (z-\lam) W(u(z),v(z)).
\ee
Similarly, $\check{H} - \lam = A_\theta A_\theta^*$ is associated with the potential
\be\label{eq:check_q}
\check{q}(x) = q(x) - 2 \frac{d}{dx} \frac{\theta'(\lam,x)}{\theta(\lam,x)}
\ee
and
\be
\check{u}(z,x) = -a_\theta u(z,x) = \frac{W_x(\theta(\lam),u(z))}{\theta(\lam,x)}
\ee
solves $\check{\tau} \check{u} = z \check{u}$ and given two solutions $u(z)$ and $v(z)$ we have
\be\label{eq:wronski2}
W(\check{u}(z),\check{v}(z)) = (z-\lam) W(u(z),v(z)).
\ee

\begin{theorem}
Assume Hypothesis~\ref{hyp:sc} holds.
Then the operator $\hat{H}= A_\phi A_\phi^* +\lam$ has an entire system of solutions
\begin{align}\label{def:hphi}
\hat{\phi}(z,x) &=\frac{1}{z-\lam}a_\phi\phi(z,x)= -\frac{W_x(\phi(\lam),\phi(z))}{(z-\lam)\phi(\lam,x)}= \frac{\int_a^x \phi(\lam,y) \phi(z,y) dy}{\phi(\lam,x)},\\\label{def:htheta}
\hat{\theta}(z,x) &=a_\phi\theta(z,x)= -\frac{W_x(\phi(\lam),\theta(z))}{\phi(\lam,x)},
\end{align}
which satisfy $W(\hat{\theta}(z),\hat{\phi}(z))=1$. In particular, $\hat{H}$ satisfies again Hypothesis~\ref{hyp:gen}.

Furthermore, the Weyl solutions of $\hat{H}$ are given by
\be
\hat{\phi}(z,x), \qquad
\hat{\psi}(z,x) = a_\phi\psi(z,x)=-\frac{W_x(\phi(\lam),\psi(z))}{\phi(\lam,x)} = \hat{\theta}(z,x) + \hat{M}(z) \hat{\phi}(z,x),
\ee
where
\be\label{eq:hatm}
\hat{M}(z) = (z-\lam) M(z)
\ee
is the singular Weyl function of $\hat{H}$ corresponding to the above system of solutions. Moreover, the associated spectral measures are related via
\be\label{eq:hatrho}
d\hat{\rho}(t) = (t-\lam) d\rho(t).
\ee
\end{theorem}

\begin{proof}
The first part is easy to check using \eqref{eq:wronski}. To see
\be\label{eq:wrphi}
W_x(\phi(\lam),\phi(z))= (z-\lam) \int_a^x \phi(\lam,y) \phi(z,y) dy
\ee
note that both sides have the same derivative and are both equal $0$ in the limit $x\to a$.

For the second part we set $\lam=0$ and fix $z\in\C\backslash\R$ without loss of generality. Moreover we will abbreviate
$\phi(\lam,x)=\phi_0(x)$, $\phi(z,x)=\phi(x)$, and $\psi(z,x)=\psi(x)$.

Now let $f$ be some absolutely continuous function with compact support in $(a,b)$ and $f' \in L^2$. Then invoking \eqref{resasa} one computes
using \eqref{eq:wronsk}, \eqref{def:hphi}, \eqref{def:htheta}, and integration by parts that
\begin{align*}
R_{\hat{H}}(z) f(x) =& \frac{1}{z} \left( A_\phi \int_a^b G(z,x,y) A_\phi^* f(y) dy - f(x) \right)\\
=& \frac{1}{z} \left( a_\phi \Big[\psi(x) \int_a^x \phi(y) (a_\phi^* f)(y) dy\Big] + a_\phi \Big[\phi(x) \int_x^b \psi(y) (a_\phi^* f)(y) dy\Big] - f(x) \right)\\
=& \frac{1}{z} \left( (a_\phi\psi)(x) \Big(\phi(x)f(x) + \int_a^x f(y)(a_\phi \phi)(y)  dy \Big) \right.\\
& \left. + (a_\phi\phi)(x) \Big( -\psi(x) f(x) + \int_x^b f(y) (a_\phi \psi)(y) dy\Big) - f(x) \right)\\
=& \frac{1}{z} \left( z\hat\psi(x) \int_a^x f(y)\hat\phi (y)  dy + z\hat\phi(x) \int_x^b f(y) \hat \psi(y) dy \right.\\
& \left. +\hat\psi(z)\phi(x)f(x)  -z\hat\phi(x)\psi(x)f(x) - f(x) \right)\\
=& \hat\psi(x) \int_a^x \hat\phi(y) f(y) dy + \hat\phi(x) \int_x^b \hat\psi(y) f(y) dy= \int_a^b \hat{G}(z,x,y) f(y) dy,
\end{align*}
which shows that
\be
\hat{G}(z,x,y) = \begin{cases}
\hat{\phi}(z,x) \hat{\psi}(z,y), & x \le y,\\
\hat{\phi}(z,y) \hat{\psi}(z,x), & x \ge y,
\end{cases}
\ee
is the Green's function of $\hat{H}$ since the set of $f$ under consideration is dense. Since the Green's function is
unique, it follows that $\hat{\phi}(z,x)$ and $\hat{\psi}(z,x)$ are the Weyl solutions of $\hat{H}$. Finally,
existence of the entire Weyl solution $\hat{\phi}(z,x)$ is equivalent to Hypothesis~\ref{hyp:gen} by Lemma~2.2 in \cite{kst}.
\end{proof}

Note
\be
\hat{\phi}(\lam,x) = -\frac{W_x(\phi(\lam),\dot{\phi}(\lam))}{\phi(\lam,x)} = \frac{\int_a^x \phi(\lam,y)^2 dy}{\phi(\lam,x)}, \qquad
\hat{\theta}(\lam,x) = \frac{1}{\phi(\lam,x)},
\ee
where the dot denotes differentiation with respect to $z$.

\begin{remark}
A few remarks are in order:
\begin{enumerate}
\item
If $\tau$ is regular at $a$ and $\theta(z,x)$, $\phi(z,x)$ are chosen to satisfy the boundary conditions
\begin{align*}
&\theta(z,a) = \cos(\alpha), \: \theta'(z,a) = -\sin(\alpha), \\
&\phi(z,a) = \sin(\alpha), \: \phi'(z,a) = \cos(\alpha),
\end{align*}
then $\hat{\tau}$ will be again regular if and only if $\sin(\alpha)\ne 0$. Moreover, in this case we have
\begin{align*}
&\hat{\theta}(z,a) = \sin(\alpha)^{-1}, \: \hat{\theta}'(z,a) = \cos(\alpha)\left(z-\lam-1 -\cot(\alpha)^2\right), \\
&\hat{\phi}(z,a) = 0, \: \hat{\phi}'(z,a) = \sin(\alpha).
\end{align*}
\item
If $\phi(\lam,x)$ is the principal positive solution near $a$ (i.e., if $H$ is the Friedrich's extension of $\tau$, see \cite{ro}, \cite{gzh}), that is,
\be
\int_a^c \frac{dx}{\phi(\lam,x)^2} = \infty,
\ee
then $\hat{\tau}$ will be limit point at $a$.
\item
Note that $\hat{H}$ has no eigenvalue at $\lam$: $\sig_p(\hat{H}) \cap \{ \lam\} = \emptyset$.
\end{enumerate}
\end{remark}

Recall that $N_\kappa^\infty$ is the subclass of $N_\kappa$
consisting of generalized Nevanlinna functions with $\kappa$ negative squares and having no nonreal poles and
the only generalized pole of nonpositive type at infinity.
In particular, $M\in N_\kappa^\infty$ admits an integral representation
\be \label{minkappa}
M(z) = (1+z^2)^k \int_\R \left(\frac{1}{t-z} - \frac{t}{1+t^2}\right) \frac{d\rho(t)}{(1+t^2)^k}
+ \sum_{j=0}^l a_j z^j,
\ee
where $\lam_0=\inf \sigma(H)$, $k\leq \kappa$, $l\le 2\kappa+1$,
\be\label{minkappa'}
a_j\in\R,  \quad\text{and}\quad
\int_\R (1+t^2)^{-k -1}d\rho(t)<\infty.
\ee
Without loss of generality we can assume that the representation \eqref{minkappa} is \emph{irreducible}, that is,
either $k=0$ or $\int_\R(1+t^2)^{-k} d\rho(t)=\infty$. For further definitions and details we refer to, e.g., \cite{krlan}, \cite[\S 2.2--2.3]{fl}, \cite[App.\ C]{kst2}.

\begin{corollary}
Assume Hypothesis~\ref{hyp:sc} holds. Assume also that the functions $M(z)$ and $\hat{M}(z)$ are connected by \eqref{eq:hatm}. Then
$M\in N_{\kappa}^\infty $ for some $\kappa\in\N_0$ if and only if
$\hat{M}\in  N_{\ti{\kappa}}^\infty$ with 
\be
\ti{\kappa}= \begin{cases}
\kappa, & \lim_{y\uparrow\infty}\frac{M(\I y)}{(\I y)^{2\kappa}}\in[0,\infty),\\
\kappa+1, & \lim_{y\uparrow\infty}\frac{M(\I y)}{(\I y)^{2\kappa}}\in[-\infty,0).
\end{cases}
\ee
\end{corollary}

\begin{proof}
The proof is similar to the proof of Lemmas~5.1 and 5.2 from \cite{kk}. Assume that $M\in N_{\kappa}^\infty $ for some $\kappa\in\N_0$. 
Using an irreducible integral representation \eqref{minkappa} we will show that
\begin{align*}
\hat{M}(z) &= (z-\lam)M(z)\\
&=(z-\lam)(1+z^2)^k \int_{\lambda_0}^{+\infty} \left(\frac{1}{t-z} - \frac{t}{1+t^2}\right) \frac{d\rho(t)}{(1+t^2)^k}+(z-\lam_0)\sum_{j=0}^l a_j z^j,
\end{align*}
where $\lam_0=\inf \sigma(H)$, lies in $N_{\tilde{\kappa}}^\infty$.  Using
\[
(z-\lam)\left(\frac{1}{t-z}-\frac{t}{1+t^2}\right)=(z^2+1)\frac{1}{t-z}\frac{t-\lam}{1+t^2}-\frac{z\lam +1}{1+t^2}
\]
we obtain
\[
\hat{M}(z) = (1+z^2)^{k+1} \int_{\lambda_0}^{+\infty} \left(\frac{1}{t-z} - \frac{t}{1+t^2}\right) \frac{d\hat{\rho}(t)}{(1+t^2)^{k+1}}+\sum_{j=0}^{l+1} \hat{a}_j z^j,
\]
where $d\hat{\rho}$ given by \eqref{eq:hatrho} is positive since $\lam\ge \lam_0$ by Hypothesis~\ref{hyp:sc}.

The converse implication $\hat{M}\in  N_{\ti{\kappa}}^\infty \Rightarrow M\in N_{\kappa}^\infty$ can be established analogously.
 
Finally, the connection between $\kappa$ and $\ti{\kappa}$ is straightforward from the following characterization of $N_\kappa^\infty$--functions:
Given a generalized Nevanlinna function $M$ in $N_\kappa^\infty$, the corresponding $\kappa$
is given by the multiplicity of the generalized pole at $\infty$ which is determined by the facts
that the following limits exist and take values as indicated:
\[
\lim_{y\uparrow \infty} -\frac{M(\I y)}{(\I y) ^{2\kappa-1}} \in (0,\infty],
\qquad
\lim_{y\uparrow \infty} \frac{M(\I y)}{(\I y) ^{2\kappa+1}} \in [0,\infty).
\]
\end{proof}

Similarly, we obtain

\begin{theorem}
Assume Hypothesis~\ref{hyp:sc}, $\tau$ is limit point at $a$, and let $\theta(\lam,x)>0$.
The operator $\check{H}= A_\theta A_\theta^* -\lam$ has an entire system of solutions
\begin{align}
\check{\phi}(z,x) &= -a_\theta\phi(z,x)=\frac{W_x(\theta(\lam),\phi(z))}{\theta(\lam,x)},\\
\check{\theta}(z,x) &=-\frac{1}{z-\lam}a_\theta\theta(z,x)=\frac{W_x(\theta(\lam),\theta(z))}{(z-\lam)\theta(\lam,x)},
\end{align}
which satisfy $W(\check{\theta}(z),\check{\phi}(z))=1$. In particular, $\check{H}$ satisfies again Hypothesis~\ref{hyp:gen}.

Furthermore, the Weyl solutions of $\check{H}$ are given by
\be
\check{\phi}(z,x), \qquad
\check{\psi}(z,x) = -a_\theta\psi(z,x)=\frac{W_x(\theta(\lam),\psi(z))}{(z-\lam)\theta(\lam,x)} = \check{\theta}(z,x) + \check{M}(z) \check{\phi}(z,x),
\ee
where
\be
\check{M}(z) = (z-\lam)^{-1} M(z)
\ee
is the singular Weyl function of $\check{H}$. The associated spectral measures are related via
\be
d\check{\rho}(t) = (t-\lam)^{-1} d\rho(t) - M(\lam) d\Theta(t-\lam),
\ee
where $\Theta(t)=0$ for $t< 0$ and $\Theta(t)=1$ for $t\ge 0$ is the usual step function.
Here $M(\lam)=\lim_{\eps\downarrow 0} M(\lam-\eps)$.
\end{theorem}

Note
\be
\check{\phi}(\lam,x) = \frac{1}{\theta(\lam,x)}, \quad
\check{\theta}(\lam,x) = \frac{W_x(\theta(\lam),\dot{\theta}(\lam))}{\theta(\lam,x)}.
\ee

\begin{remark}
Again a few remarks are in order:
\begin{enumerate}
\item
If $\theta(\lam,x)$ is the principal positive solution near $b$, that is,
\be
\int_c^b \frac{dx}{\theta(\lam,x)^2} = \infty,
\ee
then $\check{\tau}$ will be limit point at $b$. (Note that the principal positive solution near $a$ is
$\phi(\lam,x)$ since we assumed the limit point case at $a$.)
\item
Note that $\check{H}$ has an eigenvalue at $\lam$ unless $\theta(\lam,x)$ is the principal positive solution near $b$.
\item
Note that factorizing $\hat{H}$ using $\hat{\theta}(\lam,x)= \phi(\lam,x)^{-1}$ will transform $\hat{H}$ back into $H$.
In particular, $\check{\hat{\phi}}(z,x)=\phi(z,x)$ and $\check{\hat{\theta}}(z,x)=\theta(z,x)$.
\item
Clearly this procedure can be iterated and we refer (e.g.) to Appendix~A of \cite{gstdef} for the well-known formulas.
\end{enumerate}
\end{remark}

\subsection{Example: The Coulomb Hamiltonian}
We can apply the single commutation method to the Coulomb Hamiltonian
\be
H_l = -\frac{d^2}{dx^2} + \frac{l(l+1)}{x^2} - \frac{\gam}{x},\qquad l\in\N_0,
\ee
by setting
\be
A_l f = -\frac{d}{dx} + \frac{l+1}{x} + \frac{\gam}{2(l+1)},\quad
A_l^* f = \frac{d}{dx} + \frac{l+1}{x} + \frac{\gam}{2(l+1)},\\
\ee
which gives (\cite[Thm.~10.10]{tschroe})
\be
H_l = A_l^* A_l - c_l^2, \qquad H_{l+1} = A_l A_l^* - c_l^2,
\ee
where
\be
c_l = \frac{\gam}{2(l+1)}.
\ee
In particular, the singular Weyl function is given by
\be
M_l(z) = M_0(z) \prod_{k=0}^{l-1} (z-c_k^2),
\ee
where the singular Weyl function $M_0(z)$ of $H_0$ is a Herglotz--Nevanlinna function since
$H_0$ is limit circle at $a=0$.

\begin{remark}
The singular Weyl function for this case was explicitly computed in \cite[eq. (5.11)]{fl}. Moreover, the
fact that it can be reduced to the case $l=0$ via the above product formula was also first observed in
\cite[Lemma~5.1]{fl}. In the special case $\gam=0$ we obtain $M_l(z)= \sqrt{-z} z^l$ as was first observed by
\cite{DSh_00} (see also Section~5 in \cite{kst2}).
\end{remark}

\subsection{Application to perturbed Bessel operators}
Next we want to apply these results to perturbed spherical Schr\"odinger equations which have attracted
considerable interest in the past \cite{ahm}, \cite{car}, \cite{gr} \cite{kst}, \cite{kt}. In particular, we want to mention
\cite{ahm2}, where commutation techniques were used to transfer results for $l=0$ to $l\in\N_0$.

\begin{lemma}\label{lemHBes}
Fix $l\ge -\frac{1}{2}$ and $p\in[1,\infty]$. Suppose
\be\label{defHBes}
H = -\frac{d^2}{dx^2} + \frac{l(l+1)}{x^2} + q(x), \quad x\in(0,b),
\ee
where
\be
\begin{cases}
x q(x) \in L^p(0,c), & p\in(1,\infty], l\ge -\frac{1}{2},\\
x q(x) \in L^1(0,c), &  p=1, l> -\frac{1}{2},\\
x(1-\log(x)) q(x) \in L^1(0,c), & p=1, l=-\frac{1}{2},\end{cases}
\ee
for some $c\in (0,b)$. In addition, suppose $H$ is bounded from below.
If $\tau$ is limit circle at $a=0$ we impose the usual boundary condition (corresponding to the
Friedrichs extension; see also \cite{bg}, \cite{ek})
\be
\lim_{x\to0} x^l ( (l+1)f(x) - x f'(x))=0, \qquad l\in[-\frac{1}{2},\frac{1}{2}).
\ee
Then
\be\label{eq:3.35}
\hat{H}= -\frac{d^2}{dx^2} + \frac{(l+1)(l+2)}{x^2} + \hat{q}(x), \quad x \hat{q}\in L^p(0,c)
\ee
and if $l\ge\frac{1}{2}$ then
\be\label{eq:3.36}
\check{H} = -\frac{d^2}{dx^2} + \frac{(l-1)l}{x^2} + \check{q}(x), \quad x \check{q}\in L^p(0,c).
\ee
\end{lemma}

\begin{proof}
It suffices to observe that by \cite[Lem.~3.2, Cor.~3.4]{kt}
\[
\phi(\lam,x)=x^{l+1}u_1(\lam,x),
\]
where $u_1(\lam,x),\ xu_1'(\lam,x) \in W^{1,p}(0,c)$, and $\lim_{x\to 0} x u_1'(x,\lam) =0$. Noting that $\phi(\lam,x)> 0$ for all $x\in [0,b)$, by \eqref{eq:hat_q} we get
\[
x\hat{q}(x)=xq(x)-2x\frac{d}{dx}\frac{u_1'(\lam,x)}{u_1(\lam,x)} \in L^p(0,c).
\]

Similarly, to prove \eqref{eq:3.36} it suffices to note that
\[
\theta(\lam,x)=x^{-l} u_2(\lam,x),
\]
where $u_2(\lam,x),\ xu_2'(\lam,x) \in W^{1,p}(0,c)$ and $\lim_{x\to 0} x u_2'(x,\lam) =0$.
\end{proof}

Since $H$ is limit circle at $0$ for $l\in[-\frac{1}{2},\frac{1}{2})$ and thus the singular Weyl function is
a Herglotz--Nevanlinna function in this case (cf.\ \cite[App.~A]{kst2}), we obtain by
induction:

\begin{corollary}\label{cor:3.11}
Suppose
\be
H = -\frac{d^2}{dx^2} + \frac{l(l+1)}{x^2} + q(x), \qquad x q(x) \in L^p,
\ee
where $p\in[1,\infty]$ if $l+\frac{1}{2} \not\in\N_0$ and $p\in(1,\infty]$ if $l+\frac{1}{2} \in\N_0$.
Then there is a singular Weyl function of the form
\be
M(z) = (z-\lam)^{\floor{l+1/2}} M_0(z),
\ee
where $M_0(z)$ is a Herglotz--Nevanlinna function and $\lam\le\min\sig(H)$.
Here $\floor{x}= \max \{ n \in \Z | n \leq x\}$ is the usual floor function.
The corresponding spectral measure is given by
\be
d\rho(t)=(t-\lam)^{\floor{l+1/2}} d\rho_0(t),
\ee
where the measure $\rho_0$ satisfies $\int_\R d\rho_0(t)=\infty$ and $\int_\R \frac{d\rho_0(t)}{1+t^2}<\infty$.
\end{corollary}

\begin{corollary}[\cite{kst2, kt}]\label{cor:3.12}
Assume the conditions of Corollary \ref{cor:3.11}. Then there is a singular Weyl function from the class $N_\kappa^\infty$ with $\kappa= \floor{\frac{l}{2}+\frac{3}{4}}$.
\end{corollary}
\begin{proof}
The inequality $\kappa\leq \floor{\frac{l}{2}+\frac{3}{4}}$ follows after combining Corollary~\ref{cor:3.11} with Theorem~4.2 from \cite{kst2}.

It remains to show $\kappa\geq \floor{\frac{l}{2}+\frac{3}{4}}$. To this end denote by $M_0(z)$ the Weyl function of the operator $H$ with $l\in[-1/2,1/2)$.
Then, since $H$ is the Friedrichs extension of the minimal operator in this case, $M_0(z)$ satisfies (cf.\ \cite[Proposition 4]{DM91})
\be\label{eq:3.40}
M_0(z)\to -\infty, \quad \text{as}\quad z\to -\infty.
\ee
Moreover, by \cite[Cor.~A.9]{kst2} and Hypothesis \ref{hyp:sc}, $M_0(z)$ admits the following representation
\[
M_0(z)= \re(M_0(\I))+\int_\lam^{+\infty}\left(\frac{1}{t-z}-\frac{t}{1+t^2}\right)d\rho_0(t),\quad z\notin[\lam,+\infty),
\]
where the measure satisfies $\int_\lam^{+\infty} d\rho_0(t)=\infty$ and $\int_\lam^{+\infty}(1+t^2)^{-1} d\rho_0(t)<\infty$.
Hence, by \eqref{eq:3.40} we conclude that
\be\label{eq:3.41}
\int_\lam^{+\infty}\frac{d\rho_0(t)}{1+|t|}=\infty.
\ee
Combining \eqref{eq:3.41} with Corollary \ref{cor:3.11} and Theorem 4.2 from \cite{kst2}, we arrive at the desired inequality.
\end{proof}

\begin{remark}
Corollary \ref{cor:3.12} was first established by Fulton and Langer \cite{fl} in the case when the potential $q(x)$ is analytic in a neighborhood of $x=0$ (see also \cite{kl}). In the general case, it was proven in \cite{kt} (see also \cite{kst2}). Namely, by using high energy asymptotic of $\phi(z,x)$ it was shown in \cite{kst2} that $\kappa\le \ceil{\frac{l+1}{2}}$ (for the details see Section 8 in \cite{kst2}). The equality  $\kappa=\floor{\frac{l}{2}+\frac{3}{4}}$ was proven in \cite{kt} with the help of theory of super singular perturbations and detailed analysis of solutions $\phi(z,x)$ and $\theta(z,x)$.
\end{remark}

\section{Connection with the double commutation method}
\label{do_com}

In this section we want to look at the effect of the double commutation method. We refer to \cite{fg}, \cite{gt}
for further background of this method. We will use the approach from \cite{gt} with the only difference that we include the
limiting case $\gam=\infty$ (we omit the necessary minor modifications of the proofs of \cite{gt}, cf.\ \cite[Sect.~11.6]{tjac}).

Let $H$ together with a fundamental system of solutions $\phi(z,x)$, $\theta(z,x)$ as in Lemma~\ref{lem:pt} be given.

\begin{hypothesis}\label{hyp:dc}
Assume Hypothesis~\ref{hyp:gen}. Let $\gam \in (0,\infty]$ and $\lam\in\R$ such that $\phi(\lam,x)$
satisfies the boundary condition at $b$ if $\tau$ is limit circle at $b$ (i.e., $\lam$ is an eigenvalue
with eigenfunction $\phi(\lam,x)$ if $\tau$ is limit circle at $b$).
\end{hypothesis}

Introduce
\be \label{deftphigam}
\ti{\phi}_\gam(\lam,x) = \frac{\phi(\lam,x)}{\gam^{-1} + \int_a^x \phi(\lam,y)^2 dy}.
\ee
Denote by $P$ and $P_\gam$ the orthogonal projections onto the one-dimensional spaces spanned by
$\phi(\lam,x)$ and $\ti{\phi}_\gam(\lam,x)$, respectively. Here we set the projection equal to zero if
the function is not in $L^2(a,b)$. Note that $\ti{\phi}_\gam(\lam,.)\in L^2(a,b)$ if and only if $\gam<\infty$ (\cite[Lem.~2.1]{gt}).

By \cite[Lem.~2.1]{gt} the transformation
\be
U_\gam f(x) = f(x) - \ti{\phi}_\gam(\lam,x) \int_a^x \phi(\lam,y) f(y) dy
\ee
is a unitary map from $(1-P) L^2(a,b)$ onto $(1-P_\gam) L^2(a,b)$.
Moreover, by \cite[Thm.~3.2]{gt}
\be
H_\gam (1-P_\gam) = U_\gam H U_\gam^{-1} (1- P_\gam),
\ee
where the operator $H_\gam$ is associated with the potential
\be
q_\gam(x) = q(x) - 2\frac{d^2}{dx^2} \log\left( \gam^{-1} + \int_a^x \phi(\lam,y)^2 dy\right)
\ee
and boundary conditions (if any)
\be
W_a(\ti{\phi}_\gam(\lam),f) = W_b(\ti{\phi}_\gam(\lam),f)=0.
\ee
Note that for $\gam<\infty$ the operator $H_\gam$ is limit circle at $a$ if and only $H$ is and that for $\gam=\infty$ the
operator $H_\infty$ is always limit point at $b$.

\begin{theorem}
Assume Hypothesis~\ref{hyp:dc} and let $\gam<\infty$.
The operator $H_\gam$ has an entire system of solutions
\begin{align} \nn
\phi_\gam(z,x) &= \phi(z,x) - \ti{\phi}_\gam(\lam,x) \int_a^x \phi(\lam,y) \phi(z,y) dy\\\label{defphigam}
& = \phi(z,x) +\frac{1}{z-\lam} \ti{\phi}_\gam(\lam,x) W_x(\phi(\lam), \phi(z)),\\
\theta_\gam(z,x) &= \theta(z,x) +\frac{1}{z-\lam} \left( \ti{\phi}_\gam(\lam,x) W_x(\phi(\lam), \theta(z)) + \gam \phi_\gam(z,x) \right).
\end{align}
which satisfy $W(\theta_\gam(z),\phi_\gam(z)) =1$. In particular, $H_\gam$ satisfies again Hypothesis~\ref{hyp:gen}.

Furthermore, the Weyl solutions of $H_\gam$ are given by
\begin{align}\nn
\phi_\gam(z,x),\quad
\psi_\gam(z,x) &= \psi(z,x) +\frac{1}{z-\lam} \ti{\phi}_\gam(\lam,x) W_x(\phi(\lam), \psi(z))\\
& = \theta_\gam(z,x) + M_\gam(z) \phi_\gam(z,x),
\end{align}
where
\begin{align}\nn
M_\gam(z) = M(z) - \frac{\gam}{z-\lam}
\end{align}
is the singular Weyl function of $H_\gam$. The associated spectral measures are related via
\be
d\rho_\gam(x) =d\rho(x) +\gam d\Theta(x-\lam).
\ee
\end{theorem}

\begin{proof}
It is straightforward to check that $\phi_\gam(z,x)$, $\theta_\gam(z,x)$ is an entire system of solutions whose Wronskian equals
one (cf.\ \cite[(3.14) and (3.16)]{gt}). The extra multiple of $\phi_\gam(z,x)$ has been added to $\theta_\gam(z,x)$ to
remove the pole at $z=\lam$ (cf.\ \eqref{dcptlam} below).

That $\psi_\gam$ is square integrable near $b$ follows from (cf.\ \cite[(3.15)]{gt})
\be
|\psi_\gam(z,x)|^2 = |\psi(z,x)|^2 - \frac{\gam}{|z-\lam|^2} \frac{d}{dx}
\left( \frac{|W_x(\phi(\lam),\psi(z))|^2}{1 + \gam \int_a^x \phi(\lam,y)^2 dy} \right)
\ee
using the Cauchy--Schwartz inequality since
\[
W_x(\phi(\lam),\psi(z)) = W_c(\phi(\lam),\psi(z)) + (\lam-z) \int_c^x \phi(\lam,y) \psi(z,y) dy.
\]
To show that $\psi_\gam(z,x)$ satisfies the boundary condition at $b$ if $\tau_\gam$ is limit circle at $b$
we use (cf.\ \cite[(3.16)]{gt} plus \eqref{eq:wrphi})
\be
W_x(\phi_\gam(\lam),\psi_\gam(z)) = \frac{W_x(\phi(\lam),\psi(z))}{1+ \gam \int_a^x \phi(\lam,y)^2 dy}.
\ee
\end{proof}

Note
\be\label{dcptlam}
\phi_\gam(\lam,x)= \gam^{-1} \ti{\phi}_\gam(\lam,x), \qquad
\theta_\gam(\lam,x) = \theta(\lam,x) + \ti{\phi}_\gam(\lam,x) W_x(\phi(\lam), \dot{\theta}(\lam)).
\ee
\begin{remark}
Again a few remarks are in order:
\begin{enumerate}
\item
Clearly this procedure can be iterated and we refer to Section~4 of \cite{gt} for the corresponding formulas.
\item
If $\lam$ is an eigenvalue, one could even admit $\gam\in [- \|\phi(\lam)\|^{-2},\infty)$.
\item
This procedure leaves operators of the type \eqref{defHBes} invariant. In particular, it does not
change $l$.
\end{enumerate}
\end{remark}

\begin{theorem}
Assume Hypothesis~\ref{hyp:dc} and let $\gam=\infty$.
The operator $H_\infty$ has an entire system of solutions
\begin{align}\label{defphiinf}
\phi_\infty(z,x) &= \frac{1}{z-\lam} \left( \phi(z,x) - \ti{\phi}_\infty(\lam,x) \int_a^x \phi(\lam,y) \phi(z,y) dy \right)\\
\theta_\infty(z,x) &= (z-\lam) \theta(z,x) + \ti{\phi}_\infty(\lam,x) W_x(\phi(\lam), \theta(z)),
\end{align}
which satisfy $W(\theta_\infty(z),\phi_\infty(z)) =1$. In particular, $H_\infty$ satisfies again Hypothesis~\ref{hyp:gen}.

Furthermore, the Weyl solutions of $H_\infty$ are given by
\begin{align}\nn
\phi_\infty(z,x),\quad
\psi_\infty(z,x) &= (z-\lam) \psi(z,x) + \ti{\phi}_\infty(\lam,x) W_x(\phi(\lam), \psi(z))\\
& = \theta_\infty(z,x) + M_\infty(z) \phi_\infty(z,x),
\end{align}
where
\be
M_\infty(z) = (z-\lam)^2 M(z)
\ee
is the singular Weyl function of $H_\infty$. The associated spectral measures are related via
\be
d\rho_\infty(t) = (t-\lam)^2 d\rho(t).
\ee
\end{theorem}

\begin{proof}
In the limiting case $\gam \to \infty$ the definition \eqref{defphigam} from the previous theorem would give
$\phi_\infty(\lam,x)=0$ and we simply need to remove this zero. The rest follows as in the previous theorem.
\end{proof}

Note
\be
\phi_\infty(\lam,x)= \dot{\phi}(z,x) - \ti{\phi}_\infty(\lam,x) \int_a^x \phi(\lam,y) \dot{\phi}(z,y) dy, \quad
\theta_\infty(\lam,x)= -\ti{\phi}_\infty(\lam,x)
\ee

\begin{remark}
\begin{enumerate}
\item
Again this procedure can be iterated and we refer to Section~4 of \cite{gt} for the corresponding formulas.
\item
For operators of the type \eqref{defHBes} this procedure changes $l$ to $l+2$.
\end{enumerate}
\end{remark}

\section{Examples based on the generalized B\"acklund--Darboux transformation}
\label{ExpEx}

In this section we want to look at connections with the generalized B\"acklund--Darboux transformation
(GBDT) approach (see \cite{AS-GBDT0, AS-GBDT} and the references therein).
This approach contains the double commutation method as a special case as we will show below
and there are close relations with the binary Darboux transform (see, e.g., the comparative discussion in
\cite[Section 7.2]{Ci}). Here we will use the GBDT to construct an explicit example with a
generalized Weyl function which is rational with respect to $\sqrt{z}$ and which has non-real zeros.

More specific, we want to apply the GBDT to the Schr\"odinger equation 
\begin{align}&
\tau u=zu, \quad x \in (0,\infty).
\label{rE1}
\end{align}
This case was treated in Proposition~2.2 \cite{AS-PT}, see also \cite{AS-GBDT0}.

To begin with we fix an integer $n \in \N$, two $n\times n$ matrices $A$, $S(0)$ and
two vectors $\Lam_1(0),  \Lam_2(0) \in \C^n$ such that
 \be\label{rE2}
 AS(0)-S(0)A^*=
\Lambda  (0) J \Lambda (0)^*, \quad S(0)=S(0)^*, \quad \ee
where we have set
\be
\Lam(0)= \begin{bmatrix} \Lam_1(0) & \Lam_2(0) \end{bmatrix}, \qquad
J= \begin{bmatrix}
0 & 1 \\ -1 & 0
\end{bmatrix}.
\ee
Taking $\Lam(0)$ as initial condition we define the $n \times 2$ matrix function 
$\Lam(x)=\begin{bmatrix}\Lam_1(x) & \Lam_2(x) \end{bmatrix}$
as the solution of the linear system 
\be \label{E49}
\Lambda_1^{\prime  }(x)=  A  \Lambda_{2}  (x)- \Lambda_{2}  (x)q(x), \quad
\Lambda_{2}^{\prime  }(x)=  -  \Lambda_{1}(x)
\ee
and $S(x)$ via the relation
 \begin{align}\label{E49'}&
S(x)=S(0)+\int_0^x  \Lambda_{2}(t) \Lambda_{2}(t)^*dt.
\end{align}
Note that $S(x)=S(x)^*$ as well as the identity
\be \label{BD3}
AS(x)-S(x)A^*=\Lam(x)J\Lam(x)^*
\ee
which follows from \eqref{rE2}, \eqref{E49}, and \eqref{E49'}.
Furthermore, we will assume that $S(x)>0$ for $x>0$.

Given these data we can construct the Darboux matrix
using a transfer matrix function representation.  
To this end introduce
\begin{equation}  \label{BD4}
w_{A}(z,x)=I_{2}+J \Lambda (x)^{*}S(x)^{-1}
(z I_n- A)^{-1} \Lambda (x),
\end{equation}
where $I_n$ is the $n \times n$ identity matrix and the variable $x$
is added into the transfer matrix function in Lev Sakhnovich form from \cite{SaL}.
Now one can check that \eqref{BD4}
acts as a transfer matrix, that is, for any given set of linearly independent solutions $y_0$, $y_1$ of \eqref{rE1}
we obtain a set of linearly independent solutions $\wt y = \begin{bmatrix} \wt y_0 & \wt y_1 \end{bmatrix}$ of
a transformed Schr\"odinger equation $-u''+\wt q u=zu$, where the transformed
potential $\wt q$ can be expressed explicitly in terms of  $\Lam$ and $S$, by virtue of
\be\label{E48}
\wt y(z,x)= \begin{bmatrix} 1& 0\end{bmatrix}
w_{A}(z,x)w(z,x), \quad
w(z,x)=
\begin{bmatrix}
y_0(z,x) & y_1(z,x) \\ y_0'(z,x) & y_1'(z,x)
\end{bmatrix}.
\ee
In the special case $n=1$, observe that the function $\Lam_2$ is a solution of the Schr\"odinger equation
corresponding to the value $A$ of the spectral parameter. Therefore, in the case 
$n>1$, the matrix $A$ is called a generalized matrix eigenvalue. Moreover, in the
case $n=1$ the GBDT contains the double commutation method considered in Section~\ref{do_com}:

\subsection{The double commutation method as a special case of the GBDT}
\label{Comp}
Let $n=1$ (i.e., $A$ is a scalar) and set
\begin{align}\label{E50}
A=\lambda \in \C, \quad y(z, x)=
\begin{bmatrix} \phi(z, x) & \ta(z, x) \end{bmatrix},
\quad \Lam(x)=
\begin{bmatrix}
-\phi^{\prime}(\lam, x)  & \phi(\lam, x)
\end{bmatrix}.
\end{align}
Thus, for $\phi_{\gamma}$, which is given below, formulas \eqref{E48} and \eqref{E50} imply
\begin{align}\label{E51}
\phi_{\gamma}(z,x)=\ti y(z,x)\begin{bmatrix} 1 \\ 0 \end{bmatrix}=
\begin{bmatrix} 1 & 0 \end{bmatrix}
w_{A}(z,x)\begin{bmatrix}  \phi(z, x)  \\  \phi^{\prime}(z, x)  \end{bmatrix}.
\end{align}
In view of \eqref{BD4}, rewrite
\eqref{E51} in the form
\begin{align}\nn
\phi_{\gamma}(z,x)&=\phi(z, x)
+\frac{\ol \phi(\lam,x)}{S(x)(z-\lam)}\Lambda (x)
\begin{bmatrix}  \phi(z, x)  \\  \phi^{\prime}(z, x)  \end{bmatrix}
\\ \label{E52} &
=\phi(z, x)
+\frac{\ol \phi(\lam,x)}{S(x)(z-\lam)}W_x\big(\phi(\lam), \phi(z)\big).
\end{align}
By \eqref{rE2}, \eqref{E49'}, and \eqref{E50} we get
 \begin{align}\label{E53}
S(x)=(\lam - \ol \lam)^{-1}
\Lambda  (0) J \Lambda (0)^*+\int_0^x |\phi(\lam,t)|^2dt
\quad {\mathrm{for}} \quad \lam \not=\ol \lam.
\end{align}
In the case of $\lam =\ol \lam$ and real  $\phi$ treated in Section~\ref{do_com},
identity \eqref{BD3} is fulfilled automatically. Hence, noting that $S(0)\ge 0$, equality \eqref{E52} becomes
\begin{align}\label{E54}
\phi_{\gamma}(z, x) =\phi_{\gamma}(z, x)+\frac{\ti \phi_{\gamma}(\lam, x)}{z-\lam} W_x(\phi(\lam),\phi(z)),\quad \gamma:=S(0)^{-1} \in (0, \, \infty],
\end{align}
where
\begin{align} \label{E55} &
\ti \phi_{\gamma}(\lam, x)
=\frac{\phi(\lam,x)}{S(x)}, \quad S(x)=\gamma^{-1}+\int_0^x \phi(\lam,t)^2dt.
\end{align}
Hence, if $\lam =\ol \lam$ and $\phi(\lam,x)$ is real, the expressions for $\phi_\gam$ and $\ti\phi_\gam$
from above agree with \eqref{deftphigam} and \eqref{defphigam}, \eqref{defphiinf}.

\subsection{A generalized Weyl function with non-real zeros}
For the remainder of this section we consider the case
\begin{align} \label{rE3} &
q \equiv 0, \quad \Lam(0)J\Lam(0)^*=0, \quad S(0)=0.
\end{align}
Clearly, \eqref{rE2} holds for this choice of $\Lam(0)$ and $S(0)$.
For $q \equiv 0$ we have (see also \cite{gks}, where the case $q \equiv 0$,
$S(0)=I_n$ was treated):
\be  \label{E12}
\wt q(x)= 2 \Big( \big(\Lambda_{2}(x)^*S(x)^{-1}\Lambda_{2}(x)\big)^2
+ \Lambda_{1}(x)^{*}S(x)^{-1}\Lambda_{2}(x)+ \Lambda_{2}(x)^{*}
S(x)^{-1} \Lambda_{1}(x) \Big).
\ee
Moreover, for the case that $q \equiv 0$, one can set in \eqref{E48}:
\begin{equation}  \label{BD7}
w(z,x)=T(z)\begin{bmatrix}
\E^{\I\sqrt{z}x} & 0 \\ 0 & \E^{-\I\sqrt{z}x}
\end{bmatrix} T(z)^{-1},
\quad
T( z )=
\begin{bmatrix}
1 & 1  \\  \I \sqrt{z}  & -\I \sqrt{z}
\end{bmatrix}.
\end{equation}
In some places below we will need the sign in the square root fixed, hence
we chose the branch cut in $\sqrt{z}$  along the negative real axis
and assume that $\im(\sqrt{z}) >0$.
Later we will need the expression for $\ti y'$ (see \cite[formula (1.30)]{gks})
as well:
\begin{equation}  \label{BD7'}
\ti y^{\prime}(z,x )=
\begin{bmatrix} - \Lambda_{2}(x)^*S(x)^{-1}\Lambda_{2}(x) & 1 \end{bmatrix}
w_{A}(z,x )w(z,x).
\end{equation}

If $n=1$, without loss of generality we assume that
\be\label{eq:aa}
\Lam(0)=\begin{bmatrix}
v_1 & 1
\end{bmatrix}, \qquad \upsilon_1=\ol{\upsilon_1},
\ee 
and the second relation in \eqref{rE3} is immediate.

\begin{lemma}\label{Lan1}\cite{kst2}
Let $n=1$ and the parameters $A \in \C\setminus\{0\}$ and  $\upsilon_1 \in\R$ be  fixed.
Then the relation \eqref{E12}, where
 \begin{eqnarray}\nn &&
\Lam_1(x)= \upsilon_1\cos (\sqrt{A} x) +\sqrt{A}\sin (\sqrt{A} x), \\
&& \nn \Lam_2(x)= \cos (\sqrt{A} x) -\upsilon_1\sin (\sqrt{A} x)/\sqrt{A}, \quad
S(x)=\int_0^x |\Lam_2(t)|^2dt,
\end{eqnarray}
explicitly determines a singular potential $\wt q$ satisfying $\wt q(x)=
2x^{-2}(1+O(x))$ for $x \to 0$.

The corresponding entire solutions $\ti \phi(z,x)$ and $\ti \theta(z,x)$, such that 
both solutions are real-valued on $\R$,
$W(\ti \theta(z), \ti \phi(z))=1$ and $\ti \phi(z,x)$ is nonsingular
at $x=0$, are given by
\be \nn
\ti \phi(z,x)=(z-\ol A)^{-1}\ti y(z,x)
\begin{bmatrix}
1 \\ -\upsilon_1
\end{bmatrix}, \quad \ti \theta(z,x)=-(z-A)\ti y(z,x)\begin{bmatrix} 0  \\  1 \end{bmatrix}.
\ee
The singular Weyl function corresponding to this problem has the form
\begin{align}& \label{rE5}
\ti M(z)=-(z-A)(z- \ol A)\big(\I \sqrt{z}+\upsilon_1\big)^{-1}.
\end{align}
\end{lemma}

Now, we turn to the case $n=2$. To simplify calculations set
\be\label{eq_ic}
A=\begin{bmatrix}
\mu & 1 \\ 0 & \mu
\end{bmatrix}, \quad \mu \not= \ol \mu;
\quad
\Lam(0)=\begin{bmatrix} d \upsilon & \upsilon \end{bmatrix}, 
\quad  \upsilon= \begin{bmatrix} 0 \\ 1 \end{bmatrix},
\quad d=\ol d.
\ee
By \eqref{eq_ic} the second relation in \eqref{rE2} is true and we assume
also $S(0)$=0.
It is straightforward to check that \eqref{E49} (with $q \equiv 0$) and \eqref{eq_ic} hold for
\begin{align}
\label{E11}&
 \Lam_1(x)=\frac{\I}{2}\cA\left(\E^{- \om \, x}\begin{bmatrix}
c_3x+c_4 \\ c_1
\end{bmatrix}-\E^{ \om \, x}\begin{bmatrix}
c_5x- c_4  \\ c_2
\end{bmatrix}\right), \\
\label{E8}&
 \Lam_2(x)=\frac{1}{2}\left(\E^{- \om \, x}\begin{bmatrix}
c_3x+c_4 \\ c_1
\end{bmatrix}+\E^{ \om \, x}\begin{bmatrix}
c_5x-c_4  \\ c_2
\end{bmatrix}\right),
\end{align}
where $\om = \I \sqrt{\mu}$,
\begin{align} \label{E9}&
\cA =-\I \om I_2+\frac{i}{2\om}\begin{bmatrix}
0&1 \\ 0&0
\end{bmatrix}, \quad \cA^2=A,
 \quad
c_1=1+ \frac{d}{\om}, \\
\label{E10}&
c_2=1- \frac{d}{\om}, \quad
c_3=\frac{c_1 }{2  \om},  \quad
c_4=\frac{d}{2 \om^3 } , \quad c_5=-\frac{c_2 }{2  \om}.
\end{align}
The $2 \times 2$ matrix function $S$ is defined by the equality
\be\label{BD2b}
S(x)=\int_0^x  \Lambda_{2}(t) \Lambda_{2}(t)^*dt.
\ee

\begin{lemma}\label{Lan2}
Let $n=2$ and let $A$ and $\Lam(0)$ be  fixed such that \eqref{eq_ic} holds.
Then the relations \eqref{E12} and  \eqref{E11}--\eqref{BD2b} explicitly determine a singular potential $\wt q$, such that
\begin{align} \label{E22}&
\wt q(x)= \frac{12}{x^2} \big(1+O(x)\big), \quad x \geq 0, \quad x \to 0.
\end{align}
The corresponding entire real solutions $\ti \phi(z,x)$ and $\ti \theta(z,x)$, such that $\ti \phi
(z,x)$ is nonsingular
at $x=0$ and $W(\ti \theta(z),\ti \phi(z))=1$, are given by 
\begin{align} \label{E45}&
\ti \phi(z,x)=\frac{\ti y(z,x)}{(z-\ol \mu)^{2}}\begin{bmatrix}
1 \\ -d
\end{bmatrix}, \quad \ti \ta(z,x)=-(z- \mu)^{2}\ti y(z,x)\begin{bmatrix}
0 \\ 1
\end{bmatrix} ,
\end{align}
where $\ti y$ is constructed via \eqref{E48}, \eqref{BD4}, and \eqref{BD7}.
Moreover, the corresponding singular Weyl function is given by 
\begin{align}& \label{BD31}
\ti M(z)=-(z-\mu)^2(z-
 \ol \mu)^2\big(\I \sqrt{z}+d\big)^{-1}.
\end{align}
\end{lemma}

\begin{proof}
We begin by deriving some asymptotics at $x=0$. Rewrite \eqref{BD7} as
\be\label{BD7c}
w(z,x)=\cos(\sqrt{z}x)I_2 +\sin(\sqrt{z}x)\begin{bmatrix}
0 & 1/\sqrt{z}  \\  -\sqrt{z}  & 0\end{bmatrix}
\ee
to see that $w$ is an entire matrix function of $z$ (which has real-valued entries on $\R$) and the asymptotics
\be\label{E24}
w(z,x)=\big(1-\frac{z}{2}x^2\big)I_2+x
\begin{bmatrix}
0& 1 \\  -z  & 0
\end{bmatrix}+O(x^3)
\ee
hold. 
Using \eqref{E11}--\eqref{BD2b}, after some direct calculations, we get
the asymptotics of $\Lam(x)$ and $S(x)^{\pm 1}$ at $x=0$. In particular, we derive
\begin{align} \label{E20}&
\det\, S(x)= \frac{x^6}{45}\big(1+O(x)\big),  \quad S(x)^{-1}= \frac{45}{x^6} \begin{bmatrix}
x(1+O(x)) & x^3(\frac{1}{6}+O(x)) \\ x^3(\frac{1}{6}+O(x)) &x^5(\frac{1}{20}+O(x))
\end{bmatrix}.
\end{align}
Note that because of \eqref{BD2b} and the first equality in \eqref{E20} we have
$S(x)>0$ for $x>0$. Moreover, the asymptotics of $\Lam$ and formula \eqref{E20} imply
\begin{align}&\label{E28}
\Lam_2(x)^*S(x)^{-1}= - \begin{bmatrix}
15x^{-3}(1+O(x)) & \frac{3}{2}x^{-1}(1+O(x))
\end{bmatrix}.
\end{align}
Similarly, in view of \eqref{E24}, we get
\begin{align}&\label{E25}
\Lam(x)w(z,x) =\begin{bmatrix}
x-dx^2/2+O(x^3)& x^2/2 +O(x^3)\\  d- (z-\mu)x +(z-\mu)dx^2/2+O(x^3)  & 1-(z-\mu)x^2/2+O(x^3)
\end{bmatrix}.
\end{align}
Finally, by \eqref{E12},   the asymptotics of $\Lam$, and \eqref{E20}
it is easy to show that \eqref{E22} holds true. 

Now, consider
\be\label{E23}
\ti \phi_D(z,x)= \ti y(z,x)\begin{bmatrix}
1 \\ -d
\end{bmatrix}.
\ee
Clearly, $\ti \phi_D$ satisfies $\tau \ti \phi_D=z \ti \phi_D$,
where the potential is given by \eqref{E12}, since $\tau \ti y=z \ti y$.
According to \eqref{E48}, \eqref{BD4}, \eqref{eq_ic}, and \eqref{BD7c},
$\ti \phi_D$ is meromorphic in $z$ with only possible pole at $z=\mu$.
Moreover, we shall show that $\ti \phi_D(z,x)\in L^2(0,c)$. 
For that purpose rewrite \eqref{E48} as
\begin{align}&
\ti y(z,x)=
\begin{bmatrix} 1 & 0 \end{bmatrix}w(z,x)+\Lam_2(x)^*S(x)^{-1}(zI_2-A)^{-1}[\Lam_1(x)\, \, \Lam_2(x)]w(z,x).  \label{BD6}
\end{align}
It follows from \eqref{E25} and
$(zI_2-A)^{-1} =\frac{1}{(z-\mu)^2}\begin{bmatrix}
z-\mu &1\\ 0 & z-\mu
\end{bmatrix}$ that
\begin{align}&\label{E27}
(zI_2-A)^{-1}\Lam(x)w(z,x) =\frac{1}{(z-\mu)^2}\begin{bmatrix}
d +O(x^3) &1+O(x^3)\\  d(z-\mu) +O(x) & z-\mu +O(x^2)
\end{bmatrix}.
\end{align}
Relations  \eqref{E28} and  \eqref{E23}--\eqref{E27} yield $\ti \phi(z,x)\in L^2(0,c)$.

To show that $\ti \phi_D$ is an entire function,
consider again the resolvent
\begin{align}&\label{E29}
(zI_2-A)^{-1}\Lam(x)w(z,x)\begin{bmatrix}
1\\  -d
\end{bmatrix}
 =\frac{1}{(z-\mu)^2}\begin{bmatrix}
0&1\\0&0
\end{bmatrix}\Lam(x)w(\mu,x)\begin{bmatrix}
1\\  -d
\end{bmatrix}
\\ \nn &
+\frac{1}{(z-\mu)}\left(\Lam(x)w(\mu,x)+\begin{bmatrix}
0&1\\0&0
\end{bmatrix}\Lam(x)w_z(\mu,x)
\right)\begin{bmatrix}
1\\  -d
\end{bmatrix}+\Om(z,x),
\end{align}
where $w_z(z,x)=\frac{\pa}{\pa z}w(z,x)$,
and $\Om$ is an entire vector function of $z$.
In view of  \eqref{E11}--\eqref{E10} and  \eqref{BD7c},
direct calculations show that the vector coefficients at $(z-\mu)^{-2}$
and $(z-\mu)^{-1}$ on the right-hand side of \eqref{E29} equal zero,
that is, the left-hand side of \eqref{E29} is an entire vector function.
Therefore, by \eqref{E23} and \eqref{BD6}
one can see that $\ti \phi_D (z,x)$ is an entire function too.

Let us  consider another solution
\begin{align}\label{E34}&
\ti \ta_D (z,x)=\ti y(z,x)
\begin{bmatrix}
0 \\ 1 \end{bmatrix}.
\end{align}
According to \eqref{E48},  \eqref{BD7}, \eqref{BD7'}, and \eqref{E23}
we have
\begin{align}\label{E35}&
W\big(\ti \ta_D (z), \ti \phi_D(z)\big)= -\det w_A(z,x).
\end{align}
Thus, $\det w_A(z,x)$ does not depend on $x$.
Moreover, using \eqref{BD4} and \eqref{BD3} we get
\begin{align}\label{E36}&
 w_A(\ol z,x)^*J w_A(z,x)=J,
\end{align}
that is, $|\det w_A(z,x)|=1$. Taking into account the fact
that $w_A(z,x)$ is a rational function of $z$ with the only possible pole
at $z=\mu$ (of order no greater than $2$) and recalling that the relations
$\mu \not= \ol \mu$,
\begin{align}\nn&
|\det w_A(z,x)|=1, \quad \lim_{z\to \infty}\det w_A(z,x)=1
\end{align}
hold, we derive: $\det w_A(z,x)=(z-\ol \mu)^k(z-\mu)^{-k}$ 
$\,(0\leq k \leq 4)$. Further calculations show that $k=2$:
\begin{align}\label{E37}&
\det w_A(z,x)=(z-\ol \mu)^2(z-\mu)^{-2}.
\end{align}
To prove that  $\ti \phi$ and $\ti \theta$ are entire and real, rewrite \eqref{E36} as
$w_A(\ol z,x)^*=Jw_A(z,x)^{-1}J^*$. In view of   \eqref{E37} the last equality yields
\begin{align}\nn &
\begin{bmatrix}
\ol w_{11}(\ol z,x)&  \ol w_{21}(\ol z,x)
\\ \ol w_{12}(\ol z,x) &\ol w_{22}(\ol z,x) \end{bmatrix}=(z- \mu)^{2}(z-\ol \mu)^{-2}
\begin{bmatrix}
w_{11}( z,x) &   w_{21}( z,x)
\\ w_{12}( z,x)  & w_{22}( z,x) \end{bmatrix},
\end{align}
where $w_A=:\{w_{ij}\}_{i,j=1}^2$.  In other words, we get
\begin{align} \label{E44}&
(z-\ol \mu)^{-2} w_A(z,x)=(z- \mu)^{-2}
\ol w_A(\ol z,x) .
\end{align}
Recall that $w(z,x)=\ol w(\ol z,x)$, that is, $w$ is real. Thus, by \eqref{E48} and \eqref{E44} the vector function
$(z-\ol \mu)^{-2}\ti y(z,x)$ is real for $\mu \not= \ol \mu$. So, the
functions $\ti \phi$ and $\ti \ta$, which are given by \eqref{E45},  are real.
Definitions \eqref{E45}, \eqref{E23}, and \eqref{E34} yield also
\begin{align} \label{E46}&
\ti \phi(z,x)=(z-\ol \mu)^{-2}\ti \phi_D(z,x), \quad \ti \ta(z,x)=- (z-\mu)^{2}\ti \ta_D(z,x).
\end{align}

Since $\ti \phi_D$ is an entire function,
it follows from \eqref{E46} that the real function $\ti \phi$ may have only one pole
at $z=\ol \mu$ $\, (\mu \not=\ol \mu)$. Therefore, $\ti \phi$ is an entire function.
It follows from \eqref{E46} that $\ti \ta$ is an entire function too.
Finally, \eqref{E35}, \eqref{E37}, and \eqref{E46} imply that $W\big(\ti \ta(z), \ti \phi(z)\big)=1$.
Thus, the statements of our lemma regarding $\ti \phi$ and $\ti \theta$ are proved.

Now, using $\ti \phi$ and $\ti \theta$ we can construct
explicitly a singular Weyl function $\ti M$. Observe that \eqref{BD7} yields
\begin{align} \label{E47}&
w(z,x)\begin{bmatrix}
1  \\ \I \sqrt{z}
\end{bmatrix} =\begin{bmatrix}
\E^{\I\sqrt{z}x}  \\ \I\sqrt{z}\E^{\I\sqrt{z}x}
\end{bmatrix}.
\end{align}
In view of   \eqref{E8} and \eqref{BD2b} we  calculate
that
\[
\det S(x) \sim |c_2c_5|^2(4\re(\om))^{-4} \E^{4\re(\om)x}
\quad
(x \to +\infty),
\]
 and so
the transfer function $w_A(z,x)$ given by  \eqref{BD4} behaves like $O(x^4)$. Hence, it follows
from \eqref{E48} and \eqref{E47} that
\begin{align}& \label{BD30}
\widetilde{\psi}_D(z,x)= \ti y(x,z)\begin{bmatrix}
1  \\ \I \sqrt{z}
\end{bmatrix} = \begin{bmatrix}1 & 0 \end{bmatrix} w_A(z,x)\begin{bmatrix}
\E^{\I\sqrt{z}x}  \\ \I\sqrt{z}\E^{\I\sqrt{z}x}
\end{bmatrix}\in L^2(c, \infty).
\end{align}
Therefore, by \eqref{E45}, \eqref{BD31}, and \eqref{BD30} we get
\begin{align}&\nn
 \ti \psi(z,x)=\ti \theta(z,x) +\ti M(z)\ti \phi(z,x) \in L^2(c, \infty),
 \end{align}
that is, $\ti M(z)$ given by \eqref{BD31}  is a singular Weyl function of our system.
\end{proof}

\begin{remark} According to \eqref{E22}, the GBDT generated by
a $2 \times 2$ matrix $A$ of the form \eqref{eq_ic} transforms a Schr\"odinger
operator with $l=0$ into one with $l=3$.
 \end{remark}

\bigskip
\noindent
{\bf Acknowledgments.}
A.K. acknowledges the hospitality and financial support of the Erwin 
Schr\"odinger Institute and
financial support from the IRCSET PostDoctoral Fellowship Program. In addition, we are indebted to two anonymous
referees for valuable suggestions improving the presentation of the material.

\end{document}